\documentclass[11pt,bezier]{article}
\usepackage{amsmath,amssymb,amsfonts,euscript,graphicx}

\newtheorem{preproof}{{\bf \indent Proof.}}

\newenvironment{proof}[1]{\begin{preproof}{\rm
               #1}\hfill{$\Box$}}{\end{preproof}}

\newtheorem{prop}{\bf\indent Proposition}[section]
\newtheorem{defn}[prop]{\bf\indent Definition}
\newtheorem{cor}[prop]{\bf\indent Corollary}
\newtheorem{example}[prop]{\bf\indent Example}
\newtheorem{thm}[prop]{\bf\indent Theorem}
\newtheorem{rem}{\bf\indent Remark}

\title{\bf  \large $\phi$-$\delta$-$S$-primary submodule }
\author{{\normalsize  {\sc S. Najafi${}^{\mathsf{1}}$, {\sc S. Ghalandarzadeh${}^{\mathsf{1}}$}\thanks{Corresponding author: {ghalandarzadeh@kntu.ac.ir}}, {\sc A. Ranjbar Nejad Esfahani} and {\sc F. Olia}}
}\vspace{3mm}\\
{\footnotesize{${}^{\mathsf{1}}$\it  Faculty of Mathematics, K.N. Toosi
University of Technology, Tehran, Iran}}\\
}
\date{}

\begin{document}

\maketitle
\begin{abstract}
{\small  In this paper, we introduce and investigate some properties of $\phi$-$\delta$-$S$-primary submodules, which is a generalization of the $\phi$-$\delta$-primary submodules and prime submodules in general. We extend a number of main results about $\phi$-$\delta$-primary submodule  and prime submodule into the general framework. 
 Moreover, we show that there is an one-to-one correspondence between the set of all $\phi$-$\delta$-$S$-primary submodules of an $R$-module $M$ containing $K$ and the set of all $\phi$-$\delta$-$S$-primary submodules of $\frac{M}{K}$.}
\end{abstract}
\begin{center}{\section{Introduction
}}\end{center}
\label{sec1}
Throughout this paper, $R$ is a non-trivial commutative ring with an identity element, $S$ is a multiplicatively closed subset (m.c.s) of $R$ with $1\in S$ and $M$ is an unitary $R$-module.
The main purpose of this article is to introduce and investigate some properties of the concept of $\phi$-$\delta$-$S$-primary submodules, which is a generalization of the $\phi$-$\delta$-primary submodules and prime submodules (see \cite{TEKIR}, \cite{ SHARP}).
Also, this concept is a generalization of the concept of $\phi$-$\delta$-$S$-primary ideals in \cite{JABER}. We first state the necessary definitions and preliminaries.
Suppose that $\mathcal{L}(R)$ is the lattice of all ideals of $R$ and $\mathcal{L}(M)$ is the lattice of all submodules of $M$.
Let the function $ \delta : \mathcal{L}(R) \rightarrow \mathcal{L}(R) $ be an expansion fuction on $\mathcal{L}(R)$ (\cite{ZHAO}) and $ \phi : \mathcal{L}(M) \rightarrow \mathcal{L}(M) \cup \left\{ \emptyset  \right\} $ be a reduction function on $\mathcal{L}(M)$ (\cite{ZAMANI}).
For example the functions $\delta_{id}(I)=I$ and $\delta_{rad}(I)=\sqrt{I}$ are expansion functions for every $I \in \mathcal{L}( R)$ and $\phi(N)=\emptyset$ and $\phi(N)=\{0\}$  for every $N \in \mathcal{L}(M)$ are reduction functions. Let $N \in \mathcal{L}(M)$,  $J\in \mathcal{L}(R)$ and $K$ be a submodule of $M$.
The ideal quotient $(N:_R K) $ and the submodule $(N:_M J)$ are defined as follows. $\left( {N{:_R}K} \right) = \left\{ {r \in R \mid rK \subseteq N} \right\}$  and $ \left( {N{:_M}J} \right) = \left\{ {m \in M \mid mJ \subseteq N} \right\}$.
The concepts of prime ideal and prime submodule play a fundamental role in commutative algebra.
They play the same role as prime numbers in number theory.
As such, we intend to present a new generalization of them in module theory.
In  \cite{SHARP}, a proper submodule $N$ of $M$ is called prime (primary) submodule, whenever $am\in N$ for $m\in M$ and $a\in R$, then $m \in N$ or $a\in (N:_R M) (a\in \sqrt{(N:M)} )$ (2007).
In  \cite{ZAMANI}, Zamani extended the concept of prime submodules to the concept of $\phi $-prime submodules, where $\phi$ is a
reduction function (2010).
A proper submodule $N$ of $M$ is called $\phi$-prime submodule, whenever $am\in N \backslash \phi(N)$ for $m\in M$ and $a\in R$, then $m \in N$ or $a\in (N:_R M)$.
Also in \cite{TEKIR}, Ersoy extended the concept of $\phi$-prime submodules to the concept of $\phi $-$\delta$-primary submodules, where $\phi$ is a reduction fuction and $\delta$ is an expansion function (2022).
A proper submodule $N$ of $M$ is called $\phi$-$\delta$-primary submodule, whenever $am\in N \backslash \phi(N)$ for $m\in M$ and $a\in R$, then $m \in N$ or $a\in \delta (N:_R M)$.
In this paper we introduce the concept of $\phi$-$\delta$-$S$-primary submodule, where $S$ is a (m.c.s) subset of $R$ with $1\in S$. The new proposed concept is the generalization of $\phi$-$\delta$-primary submodules. 
Let $N$ be  a proper submodule of $M$, such that $\delta(N:_{R}M)\cap S=\emptyset$. 
We say $N$ is a $\phi$-$\delta$-$S$-primary submodule associated to $s\in S$, whenever $am\in N \backslash \phi(N)$ for $m\in M$ and $a\in R$, then $sm \in N$ or $sa\in \delta (N:_R M)$.  
In particular, if $\phi(N)=\emptyset$, $\delta(I)=I$ and $S=\{1\}$, then the concepts of prime submodule and $\phi$-$\delta$-$S$-primary submodule coincide.
According to Example \ref{example}, $\phi$-$\delta$-S-primary submodule will be a generalization of the $\phi$-prime, $S$-primary and $\phi$-$\delta$-primary submodules, but the converse is not true in general (Example \ref{converse}).
In a recent study \cite{JABER}, Jaber defined the concept of $\phi$-$\delta$-$S$-primary ideals. 
A proper ideal $I$, with $I\cap S=\emptyset$, is called  $\phi$-$\delta$-$S$-primary ideal of $R$ associated to $s\in S$, whenever $ab \in I \backslash \phi(I)$ then $sa\in I$ or $sb \in \delta(I)$ (2023).
Thus, the concept of $\phi$-$\delta$-$S$-primary submodules in module theory can be considered as the equivalent definition of $\phi$-$\delta$-$S$-primary ideal in ring theory.

The paper is organized as follows. In section 2, we propose some well-known properties, results and illustrative examples about $\phi$-$\delta$-$S$-primary submodules. One of the important theorems of this section is Theorem \ref{IK}, which provides an  equivalent definition to the main definition.
Moreover, we present and analyze some results to examine the behavior of $\phi$-$\delta$-$S$-primary submodules in the fractions modules, multiplication modules and finitely generated modules.

Section 3 concerns the image and the inverse image of $\phi$-$\delta$-$S$-primary submodules. We primarly show that the image and the inverse image of $\phi$-$\delta$-$S$-primary submodules are also preserved under epimorphisms (Theorem \ref{image}).
The most significant result of this section is Theorem \ref{qotient}, which demonstrates there is an one-to-one correspondence between the set of all $\phi$-$\delta$-$S$-primary submodules of $M$ containing $K$ and the set of all $\phi$-$\delta$-$S$-primary submodules of $\frac{M}{K}$.

In the last section, we investigate the $\phi$-$\delta$-$S$-primary submodules in the direct product of modules. To this end,
we present some results concerning $\phi$-$\delta$-$S$-primary submodules in the direct product of modules. 
It is noteworthy that, if $N_i$ is $\phi_i$-$\delta_i$-$S_i$-primary submodule of $M_i$, for $i=1,2$, then $ N_1\times N_2$ is not necessarily a $\phi_\times$-$\delta_\times$-$S$-primary submodule. See Example \ref{P_1} for more details. 
\vspace{1cm}
\begin{center}{\section{Properties of $\phi$-$\delta$-S-Primary submodules
}}\end{center}
Let $R$ be a non-trivial commutative ring with unity $1 \neq 0$ and $\mathcal{L}(R)$ be the lattice of all ideals of $R$. The ideal expansion function on $\mathcal{L}(R)$ is defined as follows.
\begin{defn}(See \cite{ZHAO})
	A function $ \delta : \mathcal{L}(R) \rightarrow \mathcal{L}(R) $ is called an ideal expansion of $\mathcal{L}(R)$ if whenever $I,J$ are ideals of $R$ with $ I \subseteq J $, then $ \delta(I)\subseteq\delta(J) $ and $ I\subseteq\delta(I)$.	
\end{defn}
\par In the following, we give some examples of the ideal expansion on $\mathcal{L}(R)$.
\begin{example}{\rm
	Consider  $ I \in \mathcal{L}(R)$. \\
$\noindent$			
	(1)	$\delta_{id}(I) = I$, which $\delta$ is called the identity expansion. \\
	(2) $\delta_{rad}(I) = \sqrt{I}$, which is called the radical expansion. \\
	(3) For an arbitrary fixed $J \in \mathcal{L}(R) $, consider $\delta_{res}(I) = (I:J) $. \\
	(4)	$\delta_{ann}(I) = ann(ann(I))$. \\
	(5)  For an arbitrary fixed $J \in \mathcal{L}(R) $, concider $\delta_{J}(I) = I + J$. }
\end{example}
\par Let $\delta_1, \delta_2 : \mathcal{L}(R) \rightarrow \mathcal{L}(R) $ be two ideal expansions of $R$. We denote $\delta_1 \leq \delta_2 $, if $\delta_1(I) \subseteq \delta_2(I) $ for every $I \in \mathcal{L}(R)$. 
 It is clear that $\mathcal{L}(M)$ is complete lattice with respect to set inclusion. In the following, we recall the definition of the submodule reduction function on $\mathcal{L}(M)$.
\begin{defn}(See \cite{ZAMANI})
	Let $M$ be an $R$-module and $\mathcal{L}(M)$ be the lattice of all submodules of  $M$. A function $ \phi : \mathcal{L}(M) \rightarrow \mathcal{L}(M) \cup \left\{ \emptyset  \right\} $ is called a submodule reduction of $M$ if whenever $N,K$ are submodules of $N$ with $ N \subseteq K $, then $\phi(N)\subseteq\phi(K) $ and  	$ \phi(N) \subseteq N $.		
\end{defn}
\par Now we give some examples of the submodule reduction on $\mathcal{L}(M)$.
\begin{example}{\rm
	Let $ N \in \mathcal{L}(M)$.\\
$\noindent$	
	(1)	$\phi_{\emptyset}(N) = \emptyset$ \\
	(2)	$\phi_{0}(N) = \{0\}$ \\
	(3) $\phi_{id}(N) = N $ \\
	(4) $\phi_k(N) = (N:M)^{k-1} N$ for each $k \geq 1$. \\
	(5)	$\phi_M(N) = (N:M)M$}
\end{example}
\par  Let $\phi_1, \phi_2 : \mathcal{L}(M) \rightarrow \mathcal{L}(M)\cup \left\{ \emptyset  \right\} $ be two submodule reductions of $M$. We also denote that $\phi_1 \leq \phi_2 $, if $\phi_1(N) \subseteq \phi_2(N) $ for every $N \in \mathcal{L}(M)$.
Let $M$ be a left $R$-module and $S$ be a multiplicatively closed subset (m.c.s.) of $R$. The notion of $\phi$-$\delta$-$S$-primary ideal was defined by Jaber in \cite{JABER}. We now introduce and study the concept of   $\phi$-$\delta$-$S$-primary  submodules.

\begin{defn}
	$(1)$ Let $N$ be a proper submodule of $M$ such that $ S \cap \delta(N:M)  = \emptyset $. $N$ is called $ \delta$-$S$-primary submodule of $M$ associated to $s \in S $, if whenever $ am \in N $, then $sa \in \delta(N:M) $ or $ sm \in N $ for all $ a \in R, m \in M$.

	$(2)$ Let $N$ be a proper submodule of $M$ such that $ S \cap \delta(N:M) = \emptyset $. $N$ is called $ \phi$-$\delta$-$S$-primary submodule of $M$ associated to $S \in S $, if whenever $ am \in N \backslash \phi(N) $, then $sa \in \delta(N:M) $ or $ sm \in N $ for all $ a \in R, m \in M$.
\end{defn}
\par Notice that, if there exists $s\in \delta(N:_RM)\cap S$, then $sa\in \delta(N:M)$ for every $a\in R$, which implies that $N$ is always a $\phi$-$\delta$-$S$-primary submodule. Therefore, we assume that $\delta(N:_RM)\cap S=\emptyset$ to avoid this case.
\begin{example} \label{example} {\rm Let $M$ be an $R$-module and $N$ be a proper submodule of $M$. Then the following statements are hold.\\ (1) $N$ is a prime submodule $\Leftrightarrow $ it is a $\phi_\emptyset$-$\delta_{id}$-$\{1\}$-primary submodule (\cite{SHARP}).\\ (2) $N$ is a primary submodule  $\Leftrightarrow $ it is a $\phi_\emptyset$-$\delta_{rad}$-$\{1\}$-primary submodule (\cite{SHARP}).\\ (3) $N$ is a $\phi$-prime submodule $\Leftrightarrow$ it is a $\phi$-$\delta_{id}$-$\{1\}$ primary submodule (\cite{ZAMANI}).\\ 
(4) $N$ is a $\phi$-$\delta$-primary submodule $\Leftrightarrow$ it is a $\phi$-$\delta$-$\{1\}$-primary submodule (\cite{TEKIR}).\\
(5) $N$ is a $S$-prime submodule $\Leftrightarrow$ it is a $\phi_\emptyset$-$\delta_{id}$-$S$-primary submodule (\cite{SEVIM}).
\\ (6) $N$ is a $S$-primary submodule $\Leftrightarrow$ it is a $\phi_\emptyset$-$\delta_{rad}$-$S$-primary submodule (\cite{ANSARI}).}
\end{example}
\par In the following, we give an example of $ \phi$-$\delta$-$S$-primary submodule. 
\begin{example}{\rm
Let $R = \mathbb{Z}$, $M= \mathbb{Z}_{18}$ and $N = <3>= 3\mathbb{Z}_{18}$. It is obvious that $N$ is a proper submodule of $ M$. Let $S= \{1\}$, $J= 2\mathbb{Z}$, $\delta=\delta_{rad}$ and $\phi(K) = JK$. We have $\phi \left( {\left\langle 3 \right\rangle } \right) = 6{\mathbb{Z}_{18}}$ and $\delta(N:M)=3\mathbb{Z}$.  Now suppose that $a\in R$ and $m\in M$ such that $am\in N \backslash \phi(N)$. Then $am \in \left\{ {3,9,12} \right\} $. Hence $sa\in 3\mathbb{Z}$ or $sm \in 3\mathbb{Z}_{18}$ and so $N$ is a $ \phi$-$\delta$-$S$-primary submodules.}
\end{example}
\par We next present two examples of non $ \phi$-$\delta$-$S$-primary submodules.
\begin{example}{\rm
	Let $R =\mathbb{Z}$, $M= \mathbb{Q} \times \mathbb{Q}$ be an $R$-module and $N = 2\mathbb{Z} \times 0$ be a submodule of $M$, $S= \mathbb{Z} \backslash \{0\}$. Consider $\delta (I)= \delta_{rad}(I)$, for any ideal $I$ of $R$ and $\phi(K) = \{(0, 0)\}$ for every $K \in \mathcal{L}(M)$. It is clear that $(N:M) = \{ 0\}$. Suppose that $2 \neq p \in \mathbb{Z} $ and $s \in S$ such that $gcd(p, s) = 1$. Then $p(  \tfrac{2}{p}, 0) \in (2\mathbb{Z} \times \{0\}) \backslash \phi(2\mathbb{Z} \times \{0\} )$. Since $sp \neq 0$, so $sp \notin \delta (N:M)$ and also $s(\tfrac{2}{p}, 0) \notin N$. Therefore, $N$ is not a  $\phi_0$-$\delta_{rad}$-$S$-primary submodule of $M$.}
\end{example}

\begin{example} {\rm
Let $M={\mathbb{Z}_{{p^\infty }}}$, $R=\mathbb{Z}$ and $G_t=\left\langle {\frac{1}{{{p^t}}} +\mathbb{Z} } \right\rangle $, $t \ge  0$. Also suppose that, $S = \{ {q^n} \mid n \in {\mathbb{N}_0}\} $ such that $p\neq q$. We know that ${p^k}(\frac{1}{{{p^{t + k}}}} +\mathbb{Z} ) \in {G_t}$. Since $\mathbb{Z}$ is integral domain, so $(G_t:_\mathbb{Z} {\mathbb{Z}_{{p^\infty }}})=0$ and $\delta_{rad}(G_t:_\mathbb{Z} {\mathbb{Z}_{{p^\infty }}})\cap S=\emptyset $. But there is no $s\in S$ such that $s(\frac{1}{p{t + 1}} +\mathbb{Z} ) \in {G_t}$ and $sp^{k}\in \delta_{rad}(G_t:_\mathbb{Z} {\mathbb{Z}_{{p^\infty }}})=0$. Therefore for every $t$, $G_t$ is not $\phi_0$-$\delta_{rad}$-$S$-primary submodule of $M$.}
\end{example}
  \par Suppose that $\delta_1 $,$\delta_2 $ are two ideal expansions, such that $\delta_1 \leq \delta_2 $.\\ (1) It is clear that very $ \phi$-$\delta_1$-$S$-primary submodule of $M$ is a $ \phi$-$\delta_2$-$S$-primary submodule.  \\
  (2) Also since $\delta_{id}  (I)=I\subseteq \delta(I) $, so every $\phi$-$S$-primary submodule of $M$ is a $ \phi$-$\delta$-$S$-primary submodule of $M$. 
	\par Moreover, if $\phi_1 $,$\phi_2 $ are two submodule reduction, \\
	 (1) $\phi_1 \leq \phi_2 $, then every  $ \phi_1$-$\delta$-$S$-primary submodule of $M$ is a $ \phi_2$-$\delta$-$S$-primary submodule.\\
	  (2) It is obvious that every $ \delta$-$S$-primary submodule of $M$ is a $ \phi$-$\delta$-$S$-primary submodule.
\par One can show that, every $ \phi$-$\delta$-primary submodule of $M$ is a $ \phi$-$\delta$-$S$-primary submodule of $M$ associated to $s \in S$ with $\delta \left( {N :M} \right) \cap S = \emptyset $. But in the next example we show that the converse is not true in general.
\begin{example} \label{converse} {\rm Concider $ \mathbb{Z} $-module  $M = \mathbb{Z}_{2} \times \mathbb{Z}_{3} \times \mathbb{Z}$. Let $S = \mathbb{Z} \backslash \{0\}$ and $N = (0, 0, 0)$. Suppose that $a\in R$ and $(m,n,r)\in M$, such that $a(m,n,r)\in N$. We can concider two following cases. Case1 $r\neq 0$, then $ar=0$ and since  $\mathbb{Z}$ is integral domain, so $a=0$. Hence for every $s\in S$, $as \in \delta(N:M)=0$. Case2  $r=0$, then $6(m,n,0)\in N$. Therefore $N$ is a $ \phi_{\emptyset}$-$\delta_{id}$-$S$-primary submodule of $M$ associated to $s=6$. On the otherhand $6(1,1,0)\in N\backslash \left\{ \emptyset  \right\}$, but neither $6\in \delta(N:M)$ nor $(1,1,0)\in N $. So, $N$ is not $\phi_{\emptyset}$-$\delta_{id}$-primary submodule. }
\end{example} 
\par In the following we give some properties of $ \phi$-$\delta$-$S$-primary submodules.
\begin{prop}\label{directed fam}
	Let $ \{N_i : i\in I \}  $ be a directed family of $ \phi$-$\delta$-$S$-primary submodules of $M$ associated to $s \in S$, such that $\delta \left( { \cup {N_i}:M} \right) \cap S = \emptyset $. Then the submodule $N = \cup_{i \in I} N_i$ is $ \phi$-$\delta$-$S$-primary submodule of $M$ associated to $s \in S$.			
\begin{proof}{Let $am\in N \backslash \phi(N)$, for $a \in R, m \in M$. Suppose that $sm \notin N $ for all $s \in S$. We show that $sa \in \delta(N:M)$. Since $am \notin \phi(N)$, we have $ am\notin  \phi(N_i)$ for all $i \in I$. Let  $j \in I$ such that $am \in N_j \backslash \phi(N_j)$. Since $N_j$ is a $ \phi$-$\delta$-$S$-primary submodule of $M$ associated to $s \in S$ and $sm \notin N_j$, we have $sa \in \delta(N_j:M) \subseteq \delta(N:M)$. Hence $N$ is $ \phi$-$\delta$-$S$-primary submodule of $M$ associated to $s \in S$.}
\end{proof}
\end{prop}
\par Suppose that $\delta : \mathcal{L}(R) \rightarrow \mathcal{L}(R) $ is an ideal expansion on $\mathcal{L}(R)$ and $\phi : \mathcal{L}(M) \rightarrow \mathcal{L}(M)\cup \{\emptyset\} $ is a submodule reduction on $\mathcal{L}(M)$.
We say $\delta$ and $\phi$ have intersection property if for any two ideals $I$, $J$ of $R$ and any two submodules $N$, $K$ of $M$, 
$\delta(I \cap J) = \delta(I) \cap \delta(J)$ and $\phi(N \cap K) = \phi(N) \cap \phi(K)$. 
\par In Proposition \ref{directed fam}  , we showed that for any directed familiy of $ \phi$-$\delta$-$S$-primary submodules, their union is also  $ \phi$-$\delta$-$S$-primary submodule. We also prove that the statement is true for the intersection of a finite number of $ \phi$-$\delta$-$S$-primary submodule under the suitable conditions. 
\begin{prop}
Let submodule reduction $\phi$ and  ideal expansion $\delta$ have the intersection property. Suppose that $N_1, ... , N_r$ are  $ \phi$-$\delta$-$S$-primary submodules of $M$ associated to $s \in S$. Let $\phi(N_i)= \phi(N_j)$ and $\delta(N_i:M) = \delta(N_j:M)$ for every $i,j \in  \{1,... , r\}$. Then  $N = \cap_{i \in \{1,2, ...,r\} } N_i$ is $ \phi$-$\delta$-$S$-primary submodule of $M$ associated to $s \in S$.	
\begin{proof} {The first we note that $\delta \left( { \cap {N_i}:M} \right) \cap S = \emptyset $.
 Let $am \in N \backslash \phi(N)$, where $a \in R, m \in M$. Suppose that $sa \notin \delta(N:M)$ for all $s \in S$. So $sa \notin \delta \left( { \cap {N_i}:M} \right) = \delta \left( { \cap \left( {{N_i}:M} \right)} \right) =  \cap \delta \left( {{N_i}:M} \right)$. Hence there is $j \in \{1,...,r\}$ such that $sa \notin \delta (N_j :M)$. \\ Since  $\delta(N_i :M) = \delta(N_j :M)$ for every $i,j \in \{1...,r\}$, then for every  $i \in \{1,...,r\}$ we have $sa \notin \delta(N_i :M)$. Otherhand, $\phi(N_i) = \phi(N_j)$ and $\phi$ has the property intersection, hence for every $i\in \{1,...,r\}$, $am \in {N_i} \backslash\phi \left( {{N_i}} \right)$ and $sa \notin \delta \left( {{N_i}:M} \right)$. Since all of $N_i$ are $ \phi$-$\delta$-$S$-primary submodules of $M$ associated to $s \in S$, then for every $i$, $sm \in N_i$. So $sm \in \cap N_i$ and $N$ is  $ \phi$-$\delta$-$S$-primary submodules of $M$ associated to $s \in S$.}
\end{proof}
\end{prop}
\par In the following example, we illustrate that the condition $\delta(N_i:M) = \delta(N_j:M)$ in previous proposition, is not superfluous.  
 \begin{example}{\rm
 Concider $\mathbb{Z}$-module $M = 3\mathbb{Z}$. Let $S=\{ 1\} $ be a (m.c.s) subset and $N_1 = 6\mathbb{Z}$ and $N_2 = 15\mathbb{Z}$. Then $\delta_{rad} (N_1:M)=2\mathbb{Z}\neq 5\mathbb{Z}=\delta_{rad} (N_2:M)$ and $\delta_{rad}(N_i:M)\cap S=\emptyset$. 
 If $am\in 6\mathbb{Z}$, for $a\in \mathbb{Z}$ and $m\in 3\mathbb{Z}$, then $2\mid am$ and $3\mid am$. Hence $(2\mid a \vee  2\mid m)$ and $(3\mid a \vee 3\mid m)$. Considering all cases we conclude that, $sm \in 6\mathbb{Z}$ or $sa\in 2\mathbb{Z}$, for every $s \in S$. So $N_1$ and similarly $N_2$ are 
 $ \phi_{\emptyset}$-$\delta_{rad}$-$S$-primary submodules of $M$. Note further that, $90=6\times 15 \in N_1\cap N_2=30\mathbb{Z}$, but $6 \notin 10\mathbb{Z}=\delta_{rad}(N_1\cap N_2)$ and also $15\notin 30\mathbb{Z}$. So $N_1 \cap N_2 $ is not $ \phi_{\emptyset}$-$\delta_{rad}$-$S$-primary submodule of $M$. }
 \end{example}

\begin{prop}\label{(N:_M r)}
Let $N$ be a $ \phi$-$\delta$-$S$-primary submodule of $M$ associated to $s \in S$. Suppose that for every $r \in R \backslash \delta(N:_R M)$, $(\phi(N):_M r) \subseteq \phi(N:_M r)$ and $\delta((N:_M r):M) \cap S= \emptyset $. Then $ (N :_M r)$ is $ \phi$-$\delta$-$S$-primary submodule of $M$ associated to $s \in S$.
\begin{proof}{
Let $N$ be a $ \phi$-$\delta$-$S$-primary submodule of $M$ associated to $s \in S$. Suppose that $r \in R, m \in M$ and $am \in (N :_M r) \backslash \phi(N :_M r)$. Hence, $ ram \in N$. Since $(\phi(N):_M r) \subseteq \phi(N:_M r)$, $ram \notin \phi(N)$. Since $N$ is $ \phi$-$\delta$-$S$-primary submodule of $M$ associated to $s \in S$, then $sa \in \delta(N:M)$ or $srm \in N$. So $sa \in \delta((N:_M r) : M)$ or $sm \in (N:_M r)$. Therefore, $(N:_M r)$ is $ \phi$-$\delta$-$S$-primary submodule of $M$ associated to $s \in S$.}
\end{proof}
\end{prop}
\par In Proposition \ref{(N:_M r)}, if we replace $r \in R$ by ideal $I$ of $R$, then the statement is still valid. One can prove it as follows.
\begin{prop}
Let $N$ be a $ \phi$-$\delta$-$S$-primary submodule of $M$ associated to $s \in S$ and $I$ be an ideal of $R$ such that $I \nsubseteq \delta(N:M)$. Suppose that $(\phi(N):_M I) \subseteq \phi(N:_M I)$, $ (\delta(N:M): I) \subseteq \delta((N:M):I)$ and $\delta((N:_M I):M)\cap S= \emptyset$. Then $(N:_M I)$ is a $ \phi$-$\delta$-$S$-primary submodule of $M$ associated to $s \in S$.
\begin{proof}{
Let $am \in (N:_M I) \backslash \phi((N:_M I)) $. So $am \in (N:_M I) $ and $am \notin \phi((N:_M I)) $. Thus, $amI \subseteq N $ but $amI \nsubseteq \phi(N)$. Therefore, $<a><m>I \subseteq N$ and $<a><m>I \nsubseteq \phi(N) $. Since $N$ is a $ \phi$-$\delta$-$S$-primary submodule of $M$ associated to $s \in S$,  by Theorem \ref{IK}, $s<m>I \subseteq N$  or $s<a>I \subseteq \delta(N:M)$ which implies $s<m> \in (N:_M I) $ or $s<a> \subseteq (\delta(N:M): I) \subseteq \delta((N:M):I) $. Therfore  $sm \in (N:_M I)$ or $sa \in \delta((N:_MI):M) $.}
\end{proof}
\end{prop}
\par In the next theorem, we present an alternative definition of $ \phi$-$\delta$-$S$-primary submodules based on ideals of $R$ and submodules of $M$.
\begin{thm} \label{IK}
Let $M$ be an $R$-module and $N$ be a proper submodule of $M$ with $\delta(N:M)\cap S=\emptyset$.The following statements are equivalent:\\
$(1)$ $N$ is a $ \phi$-$\delta$-$S$-primary submodule of $M$ associated to $s \in S$.\\
$(2)$ For every $a \in R \backslash \delta((N:M): s)$, $(N:_M a) \subseteq (N:_M s)$ or $ (N:_M a)=(\phi(N) :_M a)$.\\
$(3)$ For every ideal $I$ of $R$ and every submodule $K$ of $M$, if $IK \subseteq N$ and $IK \nsubseteq \phi(N)$ then $sI \subseteq \delta(N:M)$ or $sK \subseteq N$.  
\begin{proof}{
(1)$\Rightarrow$(2) Let $a \in R \backslash \delta(N:sM)=\delta((N:M):s)$ and $m \in (N:_M a)$. So $am \in N$. We can concider the following cases:
(i) $am \in \phi(N)$. Then $m \in (\phi(N) :_M a)$. Thus $(N:_M a) \subseteq (\phi(N) :_M a)$. Since $(\phi(N) :_M a) \subseteq (N:_M a)$, therefore $(N:_M a) = (\phi(N) :_M a)$.
(ii)  $am \notin \phi(N)$. Then $sm \in N$ or $sa \in \delta((N:M)) $ which implies $m \in (N:_M s)$ or $sa \in \delta((N:M))$. But since $a\notin \delta(N:sM)$, so $(N:_M a) \subseteq (N:_M s)$.\\
(2)$\Rightarrow$(3) Let $I$ be an ideal of $R$ and $K$ be a submodule of $M$ such that $IK \subseteq N$. Assume that $sI \nsubseteq \delta((N:M))$ and $sK \nsubseteq N$. We show that $IK \subseteq \phi(N)$. Let $b\in I \backslash \delta((N:M):s)$. Then by (2), $(N:_M b) \subseteq (N:_M s)$ or $ (N:_M b) = (\phi(N) :_M b)$. Since $K \subseteq (N:b)$ and $K \nsubseteq (N:_M s)$, so we can conclude that $ (N:_M b) = (\phi(N) :_M b)$. Hence $bK \subseteq \phi(N)$. Now let $b \in I \cap \delta(N:M)$. Since $sI \nsubseteq \delta(N:M)$, so there is $a_0 \in I$ such that $a_0 \notin (\delta(N:M):s)$. It is easy to show that, $b+a_0 \notin (\delta(N:M):s)$. Due to (2) and similar above argumet we conclude that $bK \subseteq (b+a_0)K \subseteq \phi(N)$. Since it holds for every $b \in I$, we can conclude that $IK \subseteq \phi(N)$.\\
(3)$\Rightarrow$(1) Let $a \in R $ and $m \in M $, such that $am \in N\backslash \phi(N)$. Then $<a><m> \subseteq N$ and $<a><m> \nsubseteq \phi(N)$ implies that $s<a> \subseteq \delta((N:M))$ or 
$s<m> \subseteq N$. Thus, $sa \in \delta(N:M) $ or $sm \in N$. Therefore, $N$ is a $ \phi$-$\delta$-$S$-primary submodule of $M$ associated to $s \in S$.}
\end{proof}
\end{thm}
\par Suppose that $\delta$ is an expansion function on $\mathcal{L}(M)$ and $\delta_R$ is an expansion function on $\mathcal{L}(R)$ such that $\delta_R(N:M)=(\delta(N):M)$
\begin{thm} Let $M$ be a multiplication module and $P$ be a submodule of $M$ with $\delta(P:M)\cap S=\emptyset$, where $S$ is a (m.c.s) subset of $R$. Then the following are equivalent:\\$(1)$ $P$ is a $\phi$-$\delta_R$-$S$-primary submodule of $M$ associated to $s \in S$.\\ $(2)$  For every submodules $L$, $N$ of $M$ with $LN \subseteq P$ and $LN \not \subseteq \phi(P)$, then there is  a $s\in S$ such that $sN \subseteq \delta(P)$ or $sL\subseteq P$.
\begin{proof} {(1)$ \Rightarrow$(2) Suppose that $P$ is a $\phi$-$\delta_R$-$S$-primary submodule of $M$ associated to $s \in S$. Let $N$ and $L$ be submodule of $M$ with $LN \subseteq P $ and $LN \not \subseteq \phi(P)$. Since $M$ is a multiplication module , so there are ideals $I$ and $J$ of $R$  such that $L=IM$, $N=JM$ and $LN=IJM$. Hence $I(JM) \subseteq P $   and  $I(JM)\not \subseteq\phi(P)$. Since $P$ is a $\phi$-$\delta_R$-$S$-primary submodule of $M$ associated to $s \in S$, so $sI \subseteq \delta_R(P:M)=  (\delta(P):M)$ or $s(JM) \subseteq P$. Thus either $sL =sIM\subseteq (\delta(P):M)M=\delta(P)$ or $sN \subseteq P$. \\ (2)$\Rightarrow$(1) Let $ IK \subseteq P$ and $IK \not \subseteq \phi(P)$, which $I$ is an ideal of $R$ and $K$ is a submodule of $M$. Since $M$ is a multiplication  module, then there is an ideal $J$ of $R$ such that $K=JM$. Thus $IM.JM=IJM \subseteq P$ and $IM.JM=IJM \not \subseteq \phi(P)$. Hence either $s(IM)\subseteq \delta(P)$ or $s(JM)\subseteq P$. Therefore either $sI \subseteq (\delta(P):M)= \delta_R(P:M)$ or $sK\subseteq P$.}
\end{proof}
\end{thm}

\begin{thm} \label{S^2}
Let $N$ be a $ \phi$-$\delta$-$S$-primary submodule of $M$ associated to $s \in S$. Suppose that $a \in R \backslash (\delta(N:M): s^2)$. Then $(N:_M sa) = (\phi(N):_M sa)$ or $(N:_M sa) = (N:_M s)$.
\begin{proof}{
Let $m \in (N:_M sa) $, so $msa \in N$. Assume that $msa \in \phi(N)$. Therefore $m \in (\phi(N):_M sa) $ and so $(N:_M sa) \subseteq (\phi(N):_M sa) $. Also we know that $(\phi(N):_M sa) \subseteq (N:_M sa) $, since $\phi(N) \subseteq N$. This implies that $(N:_M sa) = (\phi(N):_M sa)$.
Now suppose that $msa \notin \phi(N)$. Hence $msa \in N \backslash \phi(N)$. Since $N$ is a $ \phi$-$\delta$-$S$-primary submodule of $M$ associated to $s \in S$ and $a \notin (\delta(N:M): s^2)$, consequently $sm \in N$. Therefore $m \in (N:_Ms)$. Hence $(N:_M sa) \subseteq (N:_M s)$. clearly $(N:_M s) \subseteq (N:_M sa)$, hence $(N:_M sa) = (N:_M s)$.}
\end{proof}
\end{thm}

\begin{prop}
	Let $N$ be a $ \phi$-$\delta$-$S$-primary submodule of $M$ associated to $s \in S$. For every  $a \in R \backslash (\sqrt{\delta(N:M)}: s)$, we have $(N:_M sa) = (\phi(N):_M sa) \cup  (N:_M s)$.	
\begin{proof}{
 First suppose that $a \in \sqrt{(\delta(N:M): s^2)} $. Then $a ^n \in (\delta(N:M): s^2)$ for some positive integer n and therefore $(as)^n \in \delta(N:M) $. Accordinaly $(as) \in \sqrt {\delta(N:M)}$ and so $a \in (\sqrt{\delta(N:M)}: s)$. Consequently $ \sqrt{(\delta(N:M): s^2)} \subseteq (\sqrt{\delta(N:M)}: s)$. Next suppose that $a \in R \backslash (\sqrt{\delta(N:M)}: s)$. Therefore $a \in R \backslash ({\delta(N:M)}: s^2)$.  Consequently by Proposition \ref{S^2}, the proof is complete.}
\end{proof}
\end{prop}
\begin{prop} 
Let $N$ be a $ \phi$-$\delta$-$S$-primary submodule of $M$ associated to $s \in S$ and $K$ be a submodule of $M$ with $K \nsubseteq N$. If $\delta(N:K) \cap S = \emptyset$ and $\phi(N\cap K)=\phi(N)$, then $N \cap K$ is a  $ \phi$-$\delta$-$S$-primary submodule of $K$ associated to $s \in S$.	
\begin{proof} {We note that $(N:M)\subseteq (N:K)= (N\cap K:K)$. 
Let $ak\in N \cap K \backslash \phi(N \cap K)$ for $a\in R$ and $k \in K$. Then $ak \in N\backslash \phi(N)$. Therefore $sk\in N$ or $sa\in \delta(N:M)$. Since $k \in K$, then $sk\in K$. Hence if $sk\in N$, then $sk \in N\cap K$. If $sa\in \delta(N:M)$, then $sa\in \delta(N\cap K:K)$. Therefore $ N \cap K$ is a $\phi$-$\delta$-$S$-primary submodule of $K$ associated to $s \in S$.}
\end{proof}
\end{prop}
\begin{prop} \label{ideal}
Let $N$ be a proper submodule of $M$, $\phi:\mathcal{L}(M) \to \mathcal{L}(M) \cup \left\{ \emptyset  \right\}$ be a reduction function on $\mathcal{L}(M)$ and $\phi_{R} :\mathcal{L}(R) \to \mathcal{L}(R) \cup \left\{ \emptyset  \right\}$ be a reduction function on $\mathcal{L}(R)$ such that ${\phi _R}(N:M) = (\phi (N):M)$. If $N$ is a $\phi$-$\delta$-$S$-primary of $M$ associated to $s\in S$, then $(N:_RM)$ is a $ \phi_R$-$\delta$-$S$-primary ideal of $R$ associated to $s\in S$.
\begin{proof}{
Let $ab\in (N:M) \backslash \phi_{R}(N:M)=(N:M) \backslash (\phi(N):M)$ for $a,b \in R$. So $abM\subseteq N$ and $abM \not\subseteq \phi(N)$. Let $Rb=I$ and $aM=K$, then $IK\subseteq N$ and $IK \not\subseteq \phi(N)$. By Theorem \ref{IK}, $sK\subseteq N$ or $sI \subseteq \delta(N:M)$. Therefpore $saM \subseteq N$ or $sb\in \delta(N:M)$. Hence $sa \in (N:M)\subseteq \delta(N:M)$ or $sb \in \delta(N:M)$. We note that since $\delta(N:M)\cap S=\emptyset $, then $(N:M) \cap S=\emptyset $.}
\end{proof}
\end{prop}
\par The next example shows that the converse of Proposition \ref{ideal} does not hold, in general. 
\begin{example} {\rm Concider $M=\mathbb{Q} \times \mathbb{Q}$ as $\mathbb{Z}\times \mathbb{Z}$-module, $N=p\mathbb{Z} \times (0)$ is a proper submodule of $M$ and $S=\{(1,1)\}$ is a (m.c.s) subset of $\mathbb{Z}\times \mathbb{Z}$. We note that $(N:_{\mathbb{Z}\times \mathbb{Z}}M)=\{(0,0) \}$. Therefore, $(N:_{\mathbb{Z}\times \mathbb{Z}}M)=0$ is a $\phi_0$-$\delta_{rad}$-$S$-primary ideal of $\mathbb{Z}$  by default. But $N$ is not a $\phi_0$-$\delta_{rad}$-$S$-primary submodule of $M$. Since $(p,1)(1,0) \in p\mathbb{Z} \times (0) \backslash \{0\} \times \{0\}$, but $s(p,1)\notin p\mathbb{Z}\times (0)$ and $s(1,0) \notin \delta_{rad}(p\mathbb{Z} \times (0) : \mathbb{Z} \times \mathbb{Z})=(0,0)$.}
\end{example}
\par In the following we show that the converse of Proposition \ref{ideal} holds for multiplication modules. We recall that an $R$-module $M$ is multiplicatin module, if for every submodule $N$ of $M$ there exists an ideal $I $ of $R$, such that $N=IM$ \cite{ALBAST}. Let $M$ be multiplication module and $N$, $K$ be two submodule of $M$. So there are ideals $I$ and $J$ of $R$ such that $N=IM$ and $K=JM$. Hence the product of $N$ and $K$ is denoted by $NK=IJM$.
\begin{prop} \label{multi} Let the situation be as descrebed in \ref{ideal}. Let $M$ be a multiplication $R$-module and $\delta(N:M)\cap S=\emptyset$. If $(N:_RM)$ is a $ \phi_R$-$\delta$-$S$-primary ideal of $R$ associated to $s\in S$, then $N$ is a $\phi$-$\delta$-$S$-primary of $M$ associated to $s\in S$
\begin{proof} {Let $IK\subseteq N$ and $IK \not\subseteq \phi(N)$ such that $I$ be an ideal of $R$ and $K$ be a proper submodule of $M$. Since $M$ is muitiplication module, so there is an ideal $J$ of $R$ such that $K=JM$. Hence $IJ \subseteq (N:M)$ and $IJ \not\subseteq (\phi(N):M)=\phi_	R(N:M)$. Then by \cite{JABER}, Theorem 2.20, $sJ\subseteq (N:M)$ or $sI\subseteq \delta(N:M)$. Hence $sJM \subseteq N$ which implice $sK \subseteq N$ or $sI \subseteq\delta(N:M)$. So by Theorem \ref{IK}, $N$ is a $\phi$-$\delta$-$S$-primary submodule of $M$.}
\end{proof}
\end{prop}

\par We know that every $ \delta$-$S$-primary submodule of $M$ is a $ \phi$-$\delta$-$S$-primary submodule, but in the next example we show that the converse is not true.
\begin{example} \label{delta}{\rm
	Consider $M = \mathbb{Z}_{12} $ as $ \mathbb{Z} $-module. Suppose that $S=\left\{ 1 \right\}$, $I = (N:_\mathbb{Z} M)$ and $\phi(N) = IN$ and $N = <4>$ as a proper submodule of $M$. We note that $S\cap \delta_{id}(N:M)=\emptyset$.
 Since $(N:_\mathbb{Z} M)=4\mathbb{Z}$ and $\phi(N) = (N:M)N=4\mathbb{Z}N=\{0,4,8\}$. So for $a\in R$ and $m\in M$ such that $am \in N\backslash\phi(N)$, then $am\in \emptyset$. Hence $N$ is a $ \phi$-$\delta_{id}$-$S$-primary submodule of $M$ by defullt. It is clear that $2.2 \in N$, but neither $2\in N$, nor $2\in \delta_{id}(N:M)=4\mathbb{Z}$, so $N$ is not a $\delta$-$S$-primary submodule of $M$}
\end{example}
\par In the following proposition, we put a condition on the $\phi$-$\delta$-$S$-primary submodule $N$ of an $R$-module $M$ to be $\delta$-$S$-primary submodule of $M$.
\begin{prop} \label{main}
Let $N$ be a  $\phi$-$\delta$-$S$-primary submodule of $M$ associated to $s \in S$. Suppose that $(N:M)N \nsubseteq \phi(N)$. Then $N$ is a $\delta$-$S$-primary submodule of $M$ associated to $s \in S$.
\begin{proof}{
Let $am \in N $ for $a \in R, m \in M$. We can consider two following cases:\\ Case i  $am \notin \phi(N)$, we have $sa \in \delta(N:M) $ or $ sm \in N $  for $s \in S$. Case ii $am \in \phi(N)$. If $aN \nsubseteq \phi(N) $, then $n \in N$ exists such that $an \notin \phi(N)$. Thus, $a(m+n) \in N\backslash \phi(N)$, which implies $sa \in \delta(N:M) $ or $ s(m+n) \in N $  for $s \in S$, which completes the proof.
Now suppose that $aN \subseteq \phi(N)$. Similarly, if $(N:M)m \nsubseteq \phi(N)$, then the proof is complete. So assume that $(N:M)m  \subseteq \phi(N) $. Since  $(N:M)N  \nsubseteq \phi(N) $, there exists $b \in (N:M)$ and $m' \in N$ such that $bm' \notin \phi(N)$, hence $(a+b)(m+m') \in N\backslash \phi(N)$. Since $N$ is a $ \phi$-$\delta$-$S$-primary submodule, either $s(a+b) \in \delta(N:M) $ or $s(m+m') \in N$, which implies $sa \in \delta(N:M) $ or $sm \in N $. Therefore, $N$ is a $\delta$-$S$-primary submodule of $M$ associated to $s \in S$.}
\end{proof}
\end{prop}
\par Also, Example \ref{delta} shows that $(N:M)N \nsubseteq \phi(N)$ in above proposition is not superfluous. In other words, Proposition \ref{main} states that if $N$ is a $ \phi$-$\delta$-$S$-primary submodule of $M$ associated to $s \in S$, but is not a $\delta$-$S$-primary submodule of $M$ associated to $s \in S$, then  $(N:M)N \subseteq \phi(N)$. For a multiplication module, Proposition\ref{main} can be expressed as follows.
\begin{prop} \label{N^2} Let $M$ be a multiplication module, $N$ be a proper submodule of $M$  and the situation be as described in Proposition \ref{ideal}.
 Let $N$ be a  $\phi$-$\delta$-$S$-primary submodule of $M$ associated to $s \in S$, but is not a $\delta$-$S$-primary submodule of $M$ associated to $s \in S$, then $N^{2}\subseteq \phi(N)$.
\begin{proof}{Due to Propositions \ref{ideal} and \ref{multi}, we have $N$ is a  $\phi$-$\delta$-$S$-primary submodule of $M$ associated to $s \in S$ if and only if $(N:M)$ is a $\phi_{R}$-$\delta$-$S$-primary ideal of $R$ associated to $s \in S$. So,  $N$ is not a $\delta$-$S$-primary submodule of $M$ associated to $s \in S$ if and only if $(N:M)$ is not a $\delta$-$S$-primary ideal of $R$ associated to $s \in S$. Therefore $(N:M)$ is a  $\phi_	R$-$\delta$-$S$-primary ideal of $R$ associated to $s \in S$, but is not a $\delta$-$S$-primary ideal of $R$ associated to $s \in S$. Hence by \cite{JABER}, Proposition 2.17. $(N:M)^2 \subseteq \phi_{R}(N:M)$. So,  $(N:M)^{2}M \subseteq \phi_{R}(N:M)M = (\phi(N):M)M=\phi(N)$ which emplies that $N^2 \subseteq \phi(N)$.}
\end{proof}
\end{prop}
\begin{cor}  Let $M$ be a multiplication module. Suppose that $N$ is a $ \phi$-$\delta$-$S$-primary submodule of $M$, but is not a $\delta$-$S$-primary submodule of $M$ associated to $s \in S$
, then $\sqrt{N}=\sqrt{\phi(N)}$.
\begin{proof} {
Suppose that $N$ is not $\delta$-$S$-primary submodule of $M$. Hence by Proposition \ref{N^2}, $N^2\subseteq \phi(N)\subseteq N$. So $\sqrt{N^2}\subseteq \sqrt{\phi(N)}\subseteq\sqrt{ N}$. On the otherhand $N=IM$ for some ideal $I$ of $R$ and by Theorem 2.12. \cite{ALBAST}, $\sqrt{N} =\sqrt{I}M$. So $\sqrt{N^2}=\sqrt{I^2}M = \sqrt{I}M =\sqrt{N}$ which implies $\sqrt{N}=\sqrt{\phi(N)}$.}
\end{proof}
\end{cor}

\begin{prop}
Let $\phi \neq \phi_{\emptyset}$ be a submodule reduction of $M$ such that $\phi(N) = \phi^2 (N) $ for each submodule $N$ of $M$. The following statements are equivalent: 

$(1)$ Every $ \phi$-$\delta$-$S$-primary submodule of $M$ is a $\delta$-primary submodule.

$(2)$ For any submodule $N$ of $M$ , $\phi(N)$ is a $\delta$-primary submodule of $M$ and every $\delta$-$S$-primary submodule of $M$ is $\delta$-primary.
\begin{proof}{
(1)$\Rightarrow$(2) Suppose that $am \in \phi(N) \backslash \phi(\phi(N))$ for $a \in R, m \in M$. Hence,  $am \in \phi(N) \backslash \phi^2 (N)=\emptyset $. Therefore $\phi(N)$ is  $ \phi$-$\delta$-$S$-primary submodule of $M$ associated to $s \in S$ by defult. So by (1) the proof is complete. First we note that every $\delta$-$S$-primary submodule of $M$ is $ \phi$-$\delta$-$S$-primary submodule of $M$, so by (1) the proof is complete. 

(2)$\Rightarrow$(1) Let $N$ be $ \phi$-$\delta$-$S$-primary submodule of $M$ associated to $s \in S$. We show that $N$ is a $\delta$-primary submodule of $M$. Let $a \in R, m \in M$ such that $am \in N$. Case i $am \notin \phi(N) $, then $am \in N \backslash \phi(N)$. Hence $sa \in \delta(N:M)$ or $sm \in N$ and so $N$ is $\delta$-$S$-primary submodule of $M$. Therefore by (2) the proof is complete. Case ii $am \in \phi(N)$. Since $\phi(N)$ is a $\delta$-primary submodule of $M$ then $m\in \phi(N)$ or $a \in \delta(\phi(N):M)$ which implies $sm \in \phi(N) \subseteq N$ or $sa \in \delta(\phi(N):M) \subseteq\delta(N:M)$. Thus, $N$ is a $\delta$-$S$-primary submodule of $M$ associated to $s \in S$. So by (2), $N$ is a $\delta$-primary submodule.}
\end{proof}
\end{prop}

\begin{thm}
Let $N$ be a submodule of $M$. Then $N$ is a $ \phi$-$\delta$-$S$-primary submodule of $M$ associated to $s \in S$ if and only if $\tfrac{N}{\phi(N)}$ is a $\phi_0$-$\delta$-$S$-primary submodule of $\tfrac{M}{\phi(N)}$.		
\begin{proof}{
Assume that $N$ be $ \phi$-$\delta$-$S$-primary submodule of $M$ associated to $s \in S$. Let $a \in R $, $m + \phi(N) \in \tfrac{M}{\phi(N)}$ and  $0 \ne a(m + \phi(N)) \in \tfrac{N}{\phi(N)}$. Hence, $am + \phi(N) \in \tfrac{N}{\phi(N)}$ and $am \notin \phi(N)$. Since $N$ is  $ \phi$-$\delta$-$S$-primary submodule of $M$ associated to $s \in S$, $sm \in N$ or $sa \in \delta(N:M)$. Suppose that $sm \in N$, so $sm + \phi(N) \in \tfrac{N}{\phi(N)}$. But if $sa \in \delta(N:M)$, since $\delta(N:M)=\delta (\tfrac{N}{\phi(N)} : \tfrac{M}{\phi(N)})$, thus $sa \in \delta (\tfrac{N}{\phi(N)} : \tfrac{M}{\phi(N)})$. Therefore, $\tfrac{N}{\phi(N)}$ is a $\phi_0$-$\delta$-$S$-primary submodule of $\tfrac{M}{\phi(N)}$.\\
Conversely, Suppose that for $a \in R$ and $m \in M$, $am \in N \backslash \phi(N)$. So $0\ne am + \phi(N) \in \tfrac{N}{\phi(N)}$.
Since $\tfrac{N}{\phi(N)}$ is a $\phi$-$\delta$-$S$-primary submodule of $\tfrac{M}{\phi(N)}$, so $s(m + \phi(N)) \in \tfrac{N}{\phi(N)} $ which implies that $sm \in N$	or $sa \in \delta (\tfrac{N}{\phi(N)} : \tfrac{M}{\phi(N)})$ which implies $sa \in \delta(N:M)$.
Therefore, $N$ is $ \phi$-$\delta$-$S$-primary submodule of $M$ associated to $s \in S$.}	
\end{proof}
\end{thm}
\par Let $S_1 \subseteq S_2$ be  multiplicative subsets of $R$. Clearly, if $N$ is a $ \phi$-$\delta$-$S_1$-primary submodule of $M$ associated to $s \in S_1$, then $N$ is a $ \phi$-$\delta$-$S_2$-primary submodule of $M$ associated to $s \in S_2$. However, as we see in the next example, the converse is not hold, in general. 
\begin{example}{\rm
 Let $R=\mathbb{Z}$, $M=\mathbb{Z}[X]$ and $N=\left\langle {4X} \right\rangle =(4X)\mathbb{Z}$. Suppose that $\phi=\phi_{\emptyset}$, $\delta=\delta_{id}$. It is shown that $(N:M)=(\left\langle {4X} \right\rangle:_\mathbb{Z} \mathbb{Z}[X])=0$.
 Let $S_1=\{1\}$ and $S_2 = \{2^n \mid n \in \mathbb{N}\}$. Suppose that $rf(X)\in (4X)\mathbb{Z}$ for some $r\in \mathbb{Z}$ and $f(X) \in \mathbb{Z}[X]=M$. Then $\sum\limits_{i = 0}^n {(r{a_i}){X^i}}  = 4mX$ for $m \in \mathbb{Z}$. So $ra_i=0$, for $i \ne 1$  and $ra_1=4m$. If $r  \ne 0$, then $a_i=0$ for $i \ne 1$. So $f(X)=a_1 X$. since $2^2 f(X) \in \left\langle {4X} \right\rangle $, Therefore $\left\langle {4X} \right\rangle $ is a $\phi$-$\delta$-$S_2$-primary submodule associated to $2^2 \in S_2$. On the otherhand $2.2X \in \left\langle {4X} \right\rangle $ but, $1.2 \notin (N:M)$ and $1.2X \notin \left\langle {4X} \right\rangle $. So $\left\langle {4X} \right\rangle $ is not a $\phi$-$\delta$-$S_1$-primary submodule of $M$.}
\end{example}
\par In the following proposition, we show that under some condition on the multiplicatively closed subsets $S_1$ and $S_2$ of $R$ (with $S_1 \subseteq S_2$), every $\phi$-$\delta$-$S_2$-primary submodule can be a $\phi$-$\delta$-$S_1$-primary submodule. 
\begin{prop}\label{S^*}
Let $S_1 \subseteq S_2$	be multiplicatively closed subsets of $R$ such that for any $s \in S_2$, there exists $t \in S_2$ with $st \in S_1$. If $N$ is a $ \phi$-$\delta$-$S_2$-primary submodule of $M$ associated to $s \in S_2$, then $N$ is a $ \phi$-$\delta$-$S_1$-primary submodule of $M$ associated to $st \in S_1$.
\begin{proof}{
Let $t \in S_2 $ such that $st \in S_1$. 
Suppose that $a \in R $ and $ m \in M$ and $am \in N \backslash \phi(N)$. Then $sm \in N$ or $sa \in \delta(N:M)$. Therefore $stm \in N$ or $sta \in \delta(N:M)$. Hence, $N$ is a $ \phi$-$\delta$-$S_1$-primary submodule of $M$ associated to $st \in S_1$.}
\end{proof}
\end{prop}

\par Let $S$ be a (m.c.s) subset of $R$  with $1 \in S$. Set
$S^* = \{x  \in R : \tfrac{x}{1}  \in U(S^{-1}R)\} $. It is clear that $S^*$ is  multiplicatively closed subset of $S^{-1}R$ which contains $S$.

\begin{prop} Let $M$ be an $R$-module and $N$ be a proper submodule of $M$. Then $N$ is $\phi$-$\delta$-$S$-primary submodule of $M$ if and only if $N$ is $ \phi$-$\delta$-$S^*$-primary submodule of $M$.
\begin{proof}{
Since $S \subseteq S^*$ and $N$ is a $ \phi$-$\delta$-$S$-primary submodule of $M$ associated to $s \in S$, so $N$ is a  $ \phi$-$\delta$-$S^*$-primary submodule of $M$ associated to $s \in S^*$.\\ 
Now let $N$ be a $ \phi$-$\delta$-$S^*$-primary submodule of $M$ associated to $s \in S^*$. Suppose that $r \in S^*$, so  there is $\frac{a}{s} \in {S^{ - 1}}R $ such that $\frac{r}{1}\frac{a}{s} = \frac{1}{1}$. Consequently there is $t\in S$ such that $tra=ts\in S$. We now have $ta \in s^*$, because  $\frac{ta}{1}\frac{r}{st} = 1$. \\
On the otherhand we note that $S^* \cap \delta(N:M) = \emptyset$. Suppose that $r \in S^* \cap \delta(N:M) $, then there are $a \in R$ and $s \in S$ such that $\tfrac{r}{1} \tfrac{a}{s} = 1$. Then there is $t \in S$, where $tra = ts $ so $tra \in \delta(N:M) \cap S$ which is a contradiction. Therefore proof is complete by Theorem \ref{S^*}.}
\end{proof}
\end{prop}

Let $M$ be an $R$-module and $S \subseteq R $ be a multiplicatively closed subset of $R$ and suppose that $S^{-1}M$ is fraction module of $M$. Let $N$ be a submodule of $M$ with $\delta \left( {N:M} \right)\cap S  = \emptyset $. 
Suppose that, $\phi_S: \mathcal{L}(S^{-1}M) \rightarrow \mathcal{L}(S^{-1}M)\cup \left\{ \emptyset  \right\}$ be a submodule reduction on $\mathcal{L}(S^{-1}M)$ such that  $S^{-1}\phi(N) = \phi_S(S^{-1}N)$ for every $N \in \mathcal{L}(M)$ and $ \delta_S : \mathcal{L}(S^{-1}R) \rightarrow \mathcal{L}(S^{-1}R) $ be an ideal expansion on $S^{-1}R$ such that $\delta_S(S^{-1}I) = S^{-1}\delta(I)$ for every $I \in \mathcal{L}(R)$. 

\begin{prop}\label{S^{-1}M}
Let the situation be as desceribed above.
Suppose that $N$ be a $ \phi$-$\delta$-$S$-primary submodule of finitely generated $R$-module $M$ associated to $s \in S$. Then $S^{-1}N$ is a $ \phi_S$-$\delta_S$-$S$-primary submodule of $S^{-1}M$ associated to $s \in S$. 
\begin{proof} {Let $s = \frac{s}{1} \in {\delta _S}({S^{ - 1}}N:{S^{ - 1}}M) \cap S$. Therefore $s\in \delta_{S}(S^{-1}(N:M))={S^{ - 1}}(\delta (N:M))$. So $ \frac{s}{1}=\frac{a}{t}$ for some $a\in \delta(N:M)$ and $t \in S$. Therefore $ust=ua$ for some $u\in S$. Hence $ua \in S\cap \delta(N:M)$, which is a contradiction. So $ {\delta _S}({S^{ - 1}}N:{S^{ - 1}}M) \cap S=\emptyset$. Now let $t_1, t_2 \in S$, $a\in R$ and $m\in M$ such that $\frac{a}{t_1}\frac{m}{t_2}\in S^{-1}N\backslash \phi_{S}(S^{-1}N)$. So there are $n \in N$ and $s'\in S$ such that $\frac{am}{t_1t_2}=\frac{n}{s'}$. Therefore $uams'=ut_1t_2n \in N$ for some $u\in S$. On the otherhand it is clear that $uams' \notin \phi(N)$. But since $N$ is a $ \phi$-$\delta$-$S$-primary submodule of $M$, hence $sa\in \delta(N:M)$ and so $\frac{sa}{t_1}\in S^{-1}\delta(N:M)=\delta_S(S^{-1}(N:M)) \subseteq \delta_{S}(S^{-1}N:S^{-1}M)$ or $sums' \in N$, which implice that $\frac{sums'}{us't_2}=\frac{sm}{t_2} \in N$. Therefore, $S^{-1}N$ is a $ \phi_S$-$\delta_S$-$S$-primary submodule of $S^{-1}M$ associated to $s \in S$. }
\end{proof}
\end{prop}

\begin{rem} \label{AS}
Suppose that  $N$ be a proper submodule of $R$-module $M$, then\\
$(1)$ $N=M$ if and only if $(N:M)=R$.\\
$(2)$  ${S^{ - 1}}N = {S^{ - 1}}M$  if and only if $({S^{ - 1}}N:{S^{ - 1}}M) = {S^{ - 1}}R$
\end{rem}

\par Let the situation be as described before Proposition \ref{S^{-1}M}. Suppose that $N$ be a proper submodule of a finitely generated $R$-module $M$ such that $\delta((N:M):a)=(\delta(N:M):a)$ and $\phi(N:_Ma)=(\phi(N):_Ma)$ for every $a \in R$. Also assume that  $S^{-1}(\delta(N:M))\cap R= \delta (S^{-1}N \cap R)$. Then the following theorem holds under two mentioned conditions. 
\begin{thm}
Let $N$ be a proper submodule of finitely generated $R$-module $M$. Let $\delta(N:M) \cap S = \emptyset$ and $\delta_S(S^{-1}(N:M)) \neq S^{-1}R  $. 
 In addition, assume that $\phi(N) = (\phi(N):_M s)$ for some $s \in S$ and $(\phi(N) :_M t) \subseteq (\phi(N) :_M s) $ for all $t \in S$. Then the following statements are equivalent:\\
$(1)$ $N$ is  $ \phi$-$\delta$-$S$-primary submodule of $M$ associated to $s \in S$. \\
$(2)$ $(N:_M s)$ is  $ \phi$-$\delta$-primary submodule of $M$. \\
$(3)$ $S^{-1}N$ is  $ \phi_S$-$\delta_S$-primary submodule of $S^{-1}N$
and for all $t \in S$ we have, $(N:_M t) \subseteq (N:_M s)$ and $((N:M):_M t) \subseteq ((N:M):_M s)$. \\
$(4)$ $S^{-1}N$ is  $ \phi_S$-$\delta_S$-primary submodule of $S^{-1}M$ and  $S^{-1}N\cap M = (N:_M s)$ and $S^{-1}(N:M) \cap R = ((N:M):_M s)$.
\begin{proof}{
(1)$\Rightarrow$(2) Since $\delta(N:M) \cap S = \emptyset$, $(N:_M s) \neq M$. Let $am \in (N:_M s) \backslash \phi(N:_M s)$ for $a \in R$ and $m \in M$. By assumption, $am \in (N:_M s) \backslash (\phi(N):_M s)$. Therefore, $sam \in N \backslash \phi(N)$ but, $N$ is  $ \phi$-$\delta$-$S$-primary submodule of $M$ associated to $s \in S$, so $s^2m \in N$ or $sa \in \delta(N:M)$. Assume that $s^2 m\in \phi(N)$, then $m \in (\phi(N) :_M s^2) \subseteq (\phi(N) :_M s)$. So $sm \in \phi(N)$ and $sam \in \phi(N)$ which is a contradiction. Therefore  $s^2m \in N\backslash\phi(N)$ or $sa \in \delta(N:M)$. Case i If $s^2m \in N\backslash\phi(N)$. Then $s^3 \in \delta(N:M)$ which is a contradiction or $sm \in N$ and $m\in (N:s)$. Case ii If $sa \in \delta(N:M)$, then $a \in (\delta(N:M):_M s)=\delta((N:M):s)= \delta((N: s):M)$. 

(2)$\Rightarrow$(3) Since $\delta (N:M) \cap S = \emptyset $, $(N:M) \cap S = \emptyset $. Hence ${S^{ - 1}}(N:M) \ne {S^{ - 1}}R$, so $({S^{ - 1}}N:{S^{ - 1}}M) \ne {S^{ - 1}}R$ since $M$ is finitely generated $R$-module. Therefore by Remark \ref{AS}, ${S^{ - 1}}N \ne {S^{ - 1}}M$. Now
let $\tfrac{a}{s_1}\tfrac{m}{s_2} \in S^{-1}N \backslash \phi_S (S^{-1}N)$ for $\tfrac{a}{s_1} \in S^{-1}R$  and $\tfrac{m}{s_2} \in S^{-1}M$.  So there is $t \in S$ such that $tam \in N \backslash \phi(N)$, which implies $tam\in(N:_M s )\backslash\phi(N:_M s)$. $(N:_M s)$ is  $ \phi$-$\delta$-primary submodule of $M$, so $m \in (N:_M s)$ which implies $\tfrac{m}{s_2} \in S^{-1}N$ or $ta \in \delta((N: s):M)$ which implies $\tfrac{a}{s_1} \in S^{-1}\delta(N:M)=\delta_{S}(S^{-1}(N:M))\subseteq \delta_{S}(S^{-1}N:S^{-1}M)$. Therefore $S^{-1}N$ is  $ \phi_S$-$\delta_S$-primary submodule of $S^{-1}M$.

Now let $t \in S$ and $m\in (N:_Mt)$. If $m\in (\phi(N):t)$, then $m\in (\phi(N):s)\subseteq (N:s)$. If  $m\notin(\phi(N):t)$, then $mt \notin \phi(N)=(\phi(N):s)=\phi(N:s)$. Also $mt\in N\subseteq(N:s)$, so $mt \in (N:s)\backslash \phi(N:s)$. By(2) $m\in (N:s)$ which completes the proof or $t\in \delta((N:s):M)=(\delta(N:M):s)$. Hence $ts\in\delta(N:M)$ which is a contradiction.

Let $r\in( (N:M):t)$, then $rM\subseteq (N:t)\subseteq(N:s)$ and $r\in ((N:M):s)$. 

(3)$\Rightarrow$(4) The first part is obvious. We note that $ (N:_M s) \subseteq S^{-1}N \cap M$. Now let  $\tfrac{m}{1} \in  S^{-1}N \cap M$, then there is $u \in S$ such that $um \in N$. By assumption $m \in (N:_M u) \subseteq (N:_M s)$. Therefore, $S^{-1}N \cap M = (N:_M s)$.\\
Let $r\in ((N:M):s)$, then $rs \in (N:M)$. Hence $r = \frac{{rs}}{s} \in {S^{ - 1}}(N:M) \cap R$. Conversely if $r = \frac{{r}}{1} \in {S^{ - 1}}(N:M) \cap R$, then there is $t\in S$ such that $tr\in (N:M)$. Hence $r\in ((N:M):t) \subseteq ((N:M):s)$.

(4)$\Rightarrow$(1) Suppose that $am \in N \backslash \phi(N)$ where $a \in R$ and $m \in m$. First we show that $\tfrac{a}{1}\tfrac{m}{1} \notin \phi_S (S^{-1}N) $. Let $\tfrac{a}{1}\tfrac{m}{1} \in \phi_S (S^{-1}N) =S^{-1}\phi(N)$. Thus $t \in S$ exists such that $tam \in \phi(N)$, so $am \in (\phi(N):_M t) \subseteq \phi(N:_M s) = \phi(N)$. Hence $am \in \phi(N)$ which is a contradiction. So $\tfrac{a}{1}\tfrac{m}{1} \in S^{-1}N \backslash \phi_S (S^{-1}N)$, but  $N_S$ is  $ \phi_S$-$\delta_S$-primary submodule of $M_S$. Then $\tfrac{m}{1} \in S^{-1}N$ or  $\tfrac{a}{1} \in  \delta_S(S^{-1}N : S^{-1}M)$. If $\tfrac{m}{1} \in S^{-1}N$, then $m \in S^{-1}N \cap M = (N:_M s)$ and $sm \in N$. If $\tfrac{a}{1} \in  \delta_S(S^{-1}N : S^{-1}M)=\delta_S(S^{-1}(N:M)) = S^{-1}\delta(N:M) $. Thus $\tfrac{a}{1} \in S^{-1}\delta(N:M) \cap R = \delta (S^{-1}(N:M) \cap R) = \delta((N:M) :_M s) = (\delta(N:M) : s)$.}
\end{proof}
\end{thm}
\section{ 
	 Image and Inverse image of $\phi$-$\delta$-$S$-Primary submodules}
\par Let $M$ and $M'$ be two $R$-module and $f:M \to M'$ be an $R$-module homomorphism. Suppose that $\phi$ is reduction function on $\mathcal{L}(M)$ and $\phi'$ is a reduction on $\mathcal{L}(M')$. Then $f$ is called $(\phi$-$\phi')$-homomorphism if  $\phi ({f^{ - 1}}(N')) = {f^{ - 1}}(\phi' (N'))$ for all $N'\in \mathcal{L}(M')$.
 \par For any  $(\phi$-$\phi')$-epimorphism we can give the following remark. 
 \begin{rem} \label{NK}
  Let $f:M \to M'$ be a non zero $(\phi$-$\phi')$-epimorphism and $N$ be a proper submodule of $M$ containing $ker(f)$, then \\$(1)$ $(N:M)=(f(N):M')$.\\
  $(2)$  $f(\phi(N))=\phi'(f(N))$.
\begin{proof} {
$(1)$ Let $r \in (N :M)$, so $rM \subseteq N$. Since f is epimorphism, $rM'=f(rM) \subseteq f(N)$ so $r \in (f(N) : M')$. Conversely let $r \in (f(N) : M' )$, so $rM' =rf(M)\subseteq f(N)$.  Then $rM \subseteq rM+ker(f) \subseteq f^{-1}f(N) = N + ker(f) = N$. Then $rM \subseteq N$ and $r \in (N:M)$ which implies $ (N:M)=(f(N) : M')$.\\
$(2)$ We know that ${f^{ - 1}}f(N) = N+ker (f)= N$, so
 $f(\phi (N))=f(\phi (f^{ - 1}f(N))) = f(f^{ - 1}(\phi' (f(N)))) =\phi' (f(N))$. }
\end{proof}
 \end{rem}
 \begin{thm} \label{image}
 Let $f:M \to M'$ be a $(\phi$-$\phi')$-epimorphism and $N$ be a proper submodule of $M$ containing $ker(f) $. Then $N$ is $\phi$-$\delta$-$S$-primary submodule of $M$ associated to $s \in S$  if and only if $f(N)$ is a $\phi'$-$\delta$-$S$-primary submodule of $M'$ associated to $s \in S$. 
\begin{proof}{
 By Remark \ref{NK}, $\delta(f(N) : M') \cap S = \emptyset$ if and only if $\delta(N:M) \cap S= \emptyset$ and $f(N)\ne M'$. 
 Let $a \in R$, $m' \in M'$ and $am' \in f(N) \backslash \phi'(f(N))=f(N) \backslash f(\phi(N)) $, (By Remark \ref{NK}). Since $f$ is epimorphism, so there is $m\in M$ such that $m'=f(m)$. Hence $f(am) \in f(N) \backslash f(\phi(N))$. Since $ker(f) \subseteq N$, so $f(am) \in f(N)$ implies that $am \in N$. On the otherhand $f(am) \notin f(\phi(N))$ which implies that $am \notin\phi(N)$. So $am \in N \backslash \phi(N)$. As $N$ is a $\phi$-$\delta$-$S$-primary submodule of $M$ associated to $s \in S$, we conclude that $sm\in N$ or  $sa \in \delta (N:M)$. Therefore $sm'=sf(m)=f(sm) \in f(N)$ or $sa\in \delta(f(N):M')$. \\ 
 Now let $f(N)$ be a $\phi'$-$\delta$-$S$-primary submodule of $M'$ associated to $s \in S$. We show that $N$ is a $\phi$-$\delta$-$S$-primary submodule of $M$ associated to $s \in S$. Let $a\in R$, $m \in M$ and $ am \in N \backslash \phi(N)$. So $af(m) \in f(N)$ and $af(m)\notin f(\phi(N))=\phi'(f(N))$. But $f(N)$ is a $\phi'$-$\delta$-$S$-primary submodule of $M'$ associated to $s \in S$. Therefore $sa\in \delta(f(N):M')=\delta(N:M)$ or $f(sm)=sf(m) \in f(N)$, consequently $sm \in f^{-1}f(N)=N$. }
\end{proof}
 \end{thm}
 
 \begin{cor}
Let $f:M \to M'$ be a $(\phi$-$\phi')$-epimprphism. Suppose that $N'$ is $\phi'$-$\delta$-$S$-primary submodule of $M'$, then $f^{-1}(N')$ is a $\phi$-$\delta$-$S$-primary submodule of $M$.
\begin{proof}{
 We know that $ker(f)=f^{-1}(\{0\}) \subseteq f^{-1}(\phi'(N'))=\phi(f^{-1}(N')) \subseteq f^{-1}(N')$. Also $f$ is an epimorphism, so $f(f^{-1}(N'))=N'$. Since $N'$ is a $\phi'$-$\delta$-$S$-primary submodule of $M'$, so by Theorem  \ref{image}, $f^{-1}(N')$ is a $\phi$-$\delta$-$S$-primary submodule of $M$. Since $(f^{-1}(N'):M) \subseteq (N':M')$, so $\delta(f^{-1}(N'):M) \cap S=\emptyset$.}
\end{proof}
 \end{cor}
\par The next theorem is  a consequence of Theorem \ref{image}.
\begin{thm} \label{qotient} Let $M$ be an $R$-module and $K$ be a proper submodule of $M$.  Suppose that $f:M \to \frac{M}{K}$ is a $(\phi$-$\phi')$-natural epimorphism.
Then, there is an one-to-one correspondence between the set of all $\phi$-$\delta$-$S$-primary submodules of $M$ containing $K$ and the set of all $\phi$-$\delta$-$S$-primary submodules of $\frac{M}{K}$.
\end{thm}

\begin{defn}
Let $N$ be  a  $ \phi$-$\delta$-$S$-primary submodule of $M$ associated to $s \in S$. Assume that $N$ is not $\delta$-$S$-primary submodule, then there are $a \in R$ and $m \in M$ such that $am \in N$ and $am \in \phi(N), sm \notin N$ and $sa \notin \delta(N :M)$. In this case we say that $(a, m)$, is called a $\phi$-$\delta$-$S$-twin zero of $N$.
\end{defn}

\begin{defn}
	Let $N$ be  a  $ \phi$-$\delta$-$S$-primary submodule of $M$. Suppose that $I$ is an ideal of $R$ and $K$ is a proper submodule of $M$. Then $N$ is called $ \phi$-$\delta$-$S$-free twin zero submodule of $M$ respect to $IK$ such that $IK \subseteq N$ and $ IK \nsubseteq \phi(N)$ and also  $(a,k)$ is not  a $ \phi$-$\delta$-$S$-twin zero of $N$ for every $a\in I$ and $k\in K$. \\ 
A submodule $N$ of $R$-module $M$ is called $ \phi$-$\delta$-$S$-free twin zero submodule of $M$, if whenever $IK \subseteq N$ with $IK \not \subseteq \phi(N)$, for some ideal $I$ of $R$ and some submodule $K$ of $M$, then $N$ is a $ \phi$-$\delta$-$S$-free twin zero with respect to $IK$.
\end{defn}

\begin{thm}
Suppose that $N$ is a $\phi$-$\delta$-$S$-primary submodule of $M$ associated to $s \in S$. 
Then $N$ is $ \phi$-$\delta$-$S$-free twin zero if and only if for every ideal $I$ of $R$ and every proper submodule $K$ of $M$ with $IK \subseteq N$ and $IK \nsubseteq \phi(N)$, either $sK \subseteq N$ or $sI \subseteq \delta(N:M)$.
\begin{proof}{
Suppose that $N$ is a  $ \phi$-$\delta$-$S$-free twin zero submodule of $M$ associated to $s \in S$. Let $IK\subseteq N$ and $IK \not\subseteq \phi(N)$, for some ideal $I$ of $R$ and proper submodule $K$ of $M$. Suppose that $sI \nsubseteq \delta(N:M)$. Then there is $a \in I$ such that $sa \in sI$ but $sa \notin \delta(N:M)$. Let $m \in K$, then $(a, m)$ is not $ \phi$-$\delta$-$S$-twin zero of $N$. If $am \notin \phi(N)$, then $am \in N \backslash \phi(N)$. Since $N$ is $ \phi$-$\delta$-$S$-primary submodule of $M$ associated to $s \in S$ and $sa \notin \delta(N:M)$, so $sm \in N$. Therefore $sK \subseteq N$, which completes the proof.  
Now let $am \in \phi(N)$. Since $sa \notin \delta(N:M)$ and $N$ is $ \phi$-$\delta$-$S$-free twin zero, thus $sm \in N$ which implies $sK \subseteq N$.
\\Conversely, let $K$ be a submodule of $M$ and $I$ be a ideal of $R$. And suppose that $IK \subseteq N$  and $IK \not \subseteq \phi(N)$ implies that either $sK \subseteq N$ or $sI \subseteq \delta(N:M)$. We prove that $N$ is a  $ \phi$-$\delta$-$S$-free twin zero with respect to $IK$. Consider $a \in I$ and $k \in K$. If $(a,k)$ is a $ \phi$-$\delta$-$S$-twin zero of $N$, then $ak \in \phi(N)$, $sk \notin N$ and $sa\notin \delta(N:M)$, which is a contradiction. So $(a, k)$ is not a $ \phi$-$\delta$-$S$-twin zero of $N$.  So $N$ is $ \phi$-$\delta$-$S$-free twin zero with respect to $IK$. Hence $N$ is $ \phi$-$\delta$-$S$-free twin zero submodule of $M$.}
\end{proof}
\end{thm}
\section{Direct product of $\phi$-$\delta$-$S$-Primary submodules}
\par Suppose that $R=R_1 \times R_2$ is a direct product of commutative rings. Let $M_i$ be an $R_i$-module and $S_i$ be a (m.c.s) subset of $R_i$, for $i=1,2$. Then $M=M_1 \times M_2$ can be viewed as an $R$-module using scalar multiplication $(r_1, r_2)(m_1 ,m_2)= (r_1m_1, r_2m_2)$ and $S=S_1 \times S_2$ can be viewed as a (m.c.s) subset of $R$. For $i=1,2$ let $\delta_i$ be an expansion function of ideals of $R_i$ and $\phi_i$ be a reduction function of submodules of $M_i$.  We define the following two functions: $\delta_{\times}(I_1 \times I_2) = \delta_1(I_1) \times \delta_2(I_2)$ and $\phi_{\times}(N_1 \times N_2) = \phi_1(N_1) \times \phi_2(N_2)$. It is easy to see that $\delta_\times$ is an ideal expansion of $R$ and $\phi_\times$ is a submodule reduction of $M$ (\cite{EROSY}).

\begin{thm} \label{M_2} Let the situation be as described.
Let $M = M_1 \times M_2$ and $N_1$ be a proper submodule of $M_1$. The following statements are equivalent:\\
$(1)$ $N_1 \times M_2$ is a $ \phi_{\times}$-$\delta_{\times}$-$S$-primary submodule of $M$ associated to $(s_1, s_2) \in S$. \\ 
$(2)$ $N_1$ is a $\delta_{1}$-$S_1$-primary submodule of $M_1$ associated to $s_1 \in S_1$ with ${\phi _2}({M_2}) \ne {M_2}$ and also $N_1 \times M_2$ is a $\delta_{\times}$-$S$-primary submodule of $M$ associated to $(s_1, s_2) \in S$.
\begin{proof} {In the first we show that $\delta_{\times}(N_{1} \times M_{2} : M_{1} \times M_{2}) \cap (S_{1 }\times S_{2}) = \emptyset$, if and only if $\delta_1 (N_1 : M_1) \cap S_1 = \emptyset$. If $r \in \delta_1 (N_1 : M_1) \cap S_1$,  then $(r, 1) \in S_1 \times S_2$ and also $(r, 1) \in \delta_1 (N_1 : M_1) \times R_2 = \delta_1 (N_1 : M_1) \times \delta_2 (M_2 : M_2)$, which is contradiction. Conversely if $(r_1, r_2)\in \delta_{\times}(N_{1} \times M_{2} : M_{1} \times M_{2}) \cap (S_{1 }\times S_{2})=\delta_1 (N_1 : M_1) \times \delta_2 (M_2 : M_2)  \cap  (S_1 \times S_2)$, then $r_1 \in \delta_1 (N_1 : M_1) \cap S_1$, which is a contradiction.\\
(1)$\Rightarrow$(2) Suppose that  $N_1 \times M_2$ is a $ \phi_{\times}$-$\delta_{\times}$-$S$-primary submodule of $M$ associated to $(s_1, s_2) \in S$. Let $a \in R_1$, $m_1 \in M_1$ and  $am_1 \in N_1$. Since $M_2\neq \phi_2(M_2)$, so there is $m_2 \in M_2 \backslash \phi(M_2)$. Hence $(a, 1)(m_1, m_2) = (am_1, m_2) \in N_1 \times M_2  \backslash \phi_{\times}(N_1 \times M_2)$. So $(s_1, s_2)(m_1, m_2) \in N_1 \times M_2$ or $(s_1, s_2)(a, 1) \in \delta_{\times}(N_1 \times M_2 : M_1 \times M_2 ) $. Since $\delta_{\times}(N_1 \times M_2 : M_1 \times M_2 ) = \delta_{1}(N_1 : M_1) \times \delta_{2}(M_2 : M_2)$,
 hence $s_1m_1 \in N_1 $ and $s_2m_2 \in M_2$ or $ s_1a \in \delta_{1}(N_1 : M_1)$ and $s_2 \in \delta_{2}(M_2 :  M_2)$. Therefore, $s_1m_1 \in N_1$ or $ s_1a \in \delta_{1}(N_1 : M_1)$, So $N_1$ is a $\delta_1$-$S_1$-primary submodule of $M$ associated to $(s_1, s_2) \in S$.
 Now suppose that  $N_1 \times M_2$ is not a $\delta_{\times}$-$S$-primary submodule of $M$ associated to $(s_1, s_2) \in S$. By Proposition \ref{main}, $(N_1 \times M_2)((N_1 \times M_2) : (M_1 \times M_2)) \subseteq \phi_{\times}(N_1 \times M_2)$, so $(N_1 \times M_2) ((N_1 : M_1) \times (M_2 : M_2) ) \subseteq \phi_{\times}(N_1 \times M_2) $. 
 Thus, $(N_1 \times M_2) ((N_1 : M_1) \times R_2 ) =(N_1(N_1:M_1)\times (M_2\times R_2))\subseteq \phi_{1}(N_1) \times \phi_{2}(M_2) $, which results $N_1(N_1 : M_1) \subseteq \phi_{1}(N_1)$ and $M_2 \subseteq \phi_{2}(M_2)$. Which is a contradiction with $M_2\neq \phi_2(M_2)$. \\
(2)$\Rightarrow$(1) It is obvious, since we know that every $\delta_{\times }$-$S$-primary submodule of $M$ is a $\phi_{\times}$-$\delta_{\times }$-$S$-primary submodule of $M$.}
\end{proof}
\end{thm}

\begin{thm} \label{N_1}
Let $M = M_1 \times M_2$ and $N_1$ be a proper submodules of $M_1$. The following statements are equivalent:\\
$(1)$ $N_1 \times M_2$ is a $ \phi_{\times}$-$\delta_{\times}$-$S$-primary submodule of $M$ associated to $(s_1, s_2) \in S$, with $ \phi_2(M_2) = M_2$. \\
$(2)$ $N_1$ is a $ \phi_{1}$-$\delta_{1}$-$S_1$-primary submodule of $M_1$ associated to $s_1 \in S_1$. 
\begin{proof} {In the first we know that $\delta_{\times}(N_1 \times M_2:M_1 \times M_2) \cap (S_1 \times S_2) = \emptyset$ if and only if $\delta_1 (N_1 : M_1) \cap S_1 = \emptyset$.\\
 (1)$\Rightarrow$(2)  Let $a_1m_1 \in N_1 \backslash \phi_{1}(N_1)$ for $a_1 \in R_1$ and $m_1 \in M_1$. Then $(a_1, 0)(m_1, 0) \in (N_1 \times M_2) \backslash \phi_{\times}(N_1 \times M_2)$. By assumption, $(s_1, s_2)(m_1, 0) \in (N_1 \times M_2)$ or $(s_1, s_2)(a_1, 0) \in  \delta_{\times}(N_1 \times M_2 : M_1 \times M_2)$, which implies $s_1m_1 \in N_1$ or $s_1a_1 \in \delta_{1}(N_1  : M_1)$. Therefore, $N_1$ is a $ \phi_{1}$-$\delta_{1}$-$S_1$-primary submodule of $M_1$ associated to $s_1 \in S_1$. \\
(2)$\Rightarrow$(1) Suppose that  $N_1$ is a $ \phi_{1}$-$\delta_{1}$-$S_1$-primary submodule of $M_1$ associated to $s_1 \in S_1$ and $ \phi_2(M_2) = M_2$. Let $(a_1, a_2)(m_1, m_2) \in N_1 \times M_2 \backslash \phi_{\times}(N_1 \times M_2) $. Hence 
$(a_1m_1, a_2m_2)  \in N_1 \times M_2 \backslash \phi_{\times}(N_1 \times M_2)$. Thus $(a_1m_1 , a_2m_2 \in N_1 \times M_2 \backslash \phi(N_1) \times M_2$. But since $a_2m_2 \in M_2 =\phi(M_2)$, so we must have $a_1m_1\in N_1$ and $a_1m_1\notin \phi(N)$. Now we have $a_1m_1 \in N_1 \backslash  \phi_1(N_1)$. Then $s_1m_1 \in N_1 $ or $ s_1a_1 \in \delta_{1}(N_1 : M_1)$ which implies $(s_1, s_2)(m_1, m_2) = (s_1m_1, s_2m_2) \in N_1 \times M_2$ or  $(s_1, s_2)(a_1, a_2)  \in \delta_{1}(N_1 : M_1) \times \delta_{2}(M_2 : M_2) = \delta_{\times}(N_1 \times M_2 : M_1 \times M_2)$. Since $ s_2a_2 \in \delta_{2}(M_2 : M_2)=R_2$. Therefore, $N_1 \times M_2$ is a $ \phi_{\times}$-$\delta_{\times}$-$S$-primary submodule of $M$ associated to $(s_1, s_2) \in S$. }
\end{proof}
\end{thm}

\begin{rem}
There are similar results for any submodule $M_1 \times N_2$ of $R$-module $M_1 \times M_2$ which $N_2$ is a submodule of $M_2$.
\end{rem}

\begin{prop}
	Let $M = M_1 \times M_2$ be an $R = R_1 \times R_2$ module and $S = S_1 \times S_2$. Suppose that $N = N_1 \times N_2$ is $\phi_{\times}$-$\delta_{\times}$-$S$-primary submodule of $M$. Then $N_1$ is $\phi_{1}$-$\delta_{1}$-$S_1$-primary submodule of $M_1$ associated to $s_1 \in S_1$ or  $N_2 $ is $\phi_{2}$-$\delta_{2}$-$S_2$-primary submodule of $M_2$ associated to $s_2 \in S_2$.
\begin{proof}{
Let $N = N_1 \times N_2$ be $\phi_{\times}$-$\delta_{\times}$-$S$-primary submodule of $M$  associated to $(s_1, s_2) \in S_1 \times S_2$. Hence $\delta_{\times}(N_1 \times N_2 : M_1 \times M_2)\cap S = \emptyset$, so $\delta_1(N_1:M_1)\cap S_1=\emptyset$ or $\delta_2(N_2:M_2)\cap S_2=\emptyset$. Assume that $\delta_2(N_2:M_2)\cap S_2=\emptyset$.
Now we show that $N_2 $ is $\phi_{2}$-$\delta_{2}$-$S_2$-primary submodule of $M_2$ associated to $s_2 \in S_2$. 
Let $r_2 m_2 \in N_2 \backslash \phi_2 (N_2)$, for $m_2 \in M_2$ and $r_2 \in R_2$. Hence, $(1, r_2)(0, m_2)=(0, r_2 m_2) \in N_1 \times N_2 \backslash \phi_{\times}(N_1 \times N_2)$. Based on assumptions, $(s_1, s_2 ) (0, m_2) \in N_1 \times N_2$ or  $(s_1, s_2 ) (1, r_2) \in\delta_{\times} (N_1 \times N_2 : M_1 \times M_2 ) = \delta_{1} (N_1 : M_1) \times \delta_{2} (N_2 : M_2)$. Then $0 \in N_1$ and $s_2 m_2 \in N_2$ or $s_1 \in \delta_{1} (N_1 : M_1) $ and $s_2 r_2 \in \delta_{2} (N_2 : M_2)$. So $s_2m_2\in N_2$ or  $s_2 r_2 \in \delta_{2} (N_2 : M_2)$. Therefore $N_2 $ is $\phi_{2}$-$\delta_{2}$-$S_2$-primary submodule of $M_2$ associated to $s_2 \in S_2$.
 Similarly, if assume that $\delta_1(N_1:M_1)\cap S_1= \emptyset$, it can proved that $N_1$ is $\phi_{1}$-$\delta_{1}$-$S_1$-primary ideal of $M_1$ associated to $s_2 \in S_2$.}
\end{proof}
\end{prop}
\begin{cor}
Suppose that $P = P_1 \times P_2$ is an ideal of $R = R_1 \times R_2$. 
 If $P = P_1 \times P_2$ is a $\phi_{\times}$-$\delta_{\times}$-$S$-primary ideal of $R$ associated to $(s_1, s_2) \in S_1 \times S_2$, then $P_1 $ is  $\phi_{1}$-$\delta_{1}$-$S_1$-primary ideal of $R_1$ associated to $s_1 \in S_1$ or $P_2 $ is $\phi_{2}$-$\delta_{2}$-$S_2$-primary ideal of $R_2$ associated to $s_2 \in S_2$.
 \end{cor}
\par It should be noted that, if $N_1$ and $N_2$ are $\phi_1$-$\delta_1$-$S_1$-primary and $\phi_2$-$\delta_2$-$S_2$-primary submodules of $M_1$ and $M_2$ respectely, then $N_1 \times N_2$ is not $\phi_{\times}$-$\delta_{\times}$-$S$-primary submodule, in general.
\begin{example} \label{P_1}{\rm
Let $M_1 \times M_2=\mathbb{Z} \times \mathbb{Z}$ be an $R_1 \times R_2=\mathbb{Z}\times \mathbb{Z}$-module. Suppose that $N_1=p\mathbb{Z}$, $N_2=q\mathbb{Z}$, such that $p$ and $q$ are two distinct prime numbers and $S_1=\{q^n \mid n\in \mathbb{N}_0  \},  S_2=\{p^n \mid n\in \mathbb{N}_0\}$. Let $am \in p\mathbb{Z} \backslash \emptyset$, so $p\mid am$. Hence $p\mid a$ or $p\mid m$ and $m\in p\mathbb{Z}$ or $a\in \delta_{rad}(p\mathbb{Z}:\mathbb{Z})=p\mathbb{Z}$. Therefore $N_1$ and similarly $N_2$ are $\phi_\emptyset$-$\delta_{rad}$-$S_i$-primary submodule of $\mathbb{Z}$. On the otherhand $(p,1)(1,q)\in p\mathbb{Z} \times q\mathbb{Z}$, but for every $(s_1,s_2) \in S_1 \times S_2$ neither $(s_1,s_2)(p,1) \in p\mathbb{Z} \times q\mathbb{Z}$, nor $(s_1,s_2)(1,q) \in \delta_\times(N_1\times N_2:M_1 \times M_2)=N_1\times N_2$. Therfore $p\mathbb{Z} \times q\mathbb{Z}$ is not a $\phi_{\times }$-$\delta_\times$-$S$-primary submodule of $M_1 \times M_2$.}
\end{example}

\end{document}